\documentclass[11pt,reqno]{article}
\usepackage{graphics}
\usepackage{epsfig}
\usepackage{amsmath}
\usepackage{epstopdf}
\usepackage{amsfonts}
\font\escript=eusm10  

\newcommand{\Cc}{\mbox{\escript C}}
\newcommand\norm[1]{\lVert#1\rVert}
\newcommand{\verti}[1]{{\left\vert\kern-0.25ex\left\vert\kern-0.25ex\left\vert#1
		\right\vert\kern-0.25ex\right\vert\kern-0.25ex\right\vert}}

\newtheorem{theorem}{Theorem}[section]
\newtheorem{lemma}[theorem]{Lemma}

\newtheorem{remark}[theorem]{Remark}

\newtheorem{example}[theorem]{Example}

\catcode`@=11 \@addtoreset{equation}{section}
\renewcommand\theequation{\thesection.\@arabic\c@equation}
\catcode`@=12

\def\eop{{\ \vrule height 7pt width 7pt depth 0pt}}
\makeatother

\begin{document}
\author{
{\sc Sandip Maji$^1$},\,\,{and }\,\,{\sc Srinivasan~Natesan$^{2,}$}\thanks{Corresponding Author, Department of Mathematics, Indian Institute of Technology Guwahati, Guwahati, Assam-781 039, {\sc India}, Email: natesan@iitg.ac.in.}\\{$^1$Department of Mathematics, Indian Institute of Technology Guwahati},\\
{ Guwahati - 781039, India,}\\{Email: s.maji@iitg.ac.in}\\{$^2$Department of Mathematics, Indian Institute of Technology Guwahati,}\\
{ Guwahati - 781039, India,}\\{Email: natesan@iitg.ac.in}}
\title{{\bf Superconvergence analysis of interior penalty discontinuous Galerkin method for a class of time-fractional diffusion problems}}

\maketitle
\thispagestyle{empty}
\bigskip
\begin{abstract}
In this study, we consider a class of non-autonomous time-fractional partial advection-diffusion-reaction (TF-ADR) equations with Caputo type fractional derivative. To obtain the numerical solution of the model problem, we apply the non-symmetric interior penalty Galerkin (NIPG) method in space on a uniform mesh and the L1-scheme in time on a graded mesh. It is demonstrated that the computed solution is discretely stable. Superconvergence of error estimates for the proposed method are obtained using the discrete energy-norm. Also, we have applied the proposed method to solve semilinear problems after linearizing by the Newton linearization process. The theoretical results are verified through numerical experiments.
\end{abstract}
\bigskip
\goodbreak\noindent
\goodbreak\noindent
\goodbreak\noindent
{\bf Author Contribution Statement:} Each author who satisfies the criteria for authorship attests that they made an adequate contribution to the work to take ownership of its content, including conception, design, analysis, writing, or manuscript change. Also, each author certifies that no publication has or will get this material or any material that is exactly the same.
\goodbreak\noindent
{\bf Conflict of Interest Statement:} There are no competing interests.
We want to reassure you that this publication has no known conflicts of interest.

\goodbreak\noindent
\section{Introduction}\label{ch5_secint}
The discontinuous Galerkin (DG) finite element methodology has become an important alternative for numerically solving advection problems for which the continuous (conforming) finite element approach fails. The fact that DG conserves mass locally is a key aspect. One can find early research on DG approaches in Reed and Hill \cite{reed1973triangular}, Arnold \cite{Arnold_SIAMJNA1982}, etc. For further details one may refer to \cite{BenLi_DGM_book,Beatrice_book_DGM}. Singh and Natesan \cite{Singh_Natesan_Calcolo18} investigated the superconvergence aspects of the NIPG approach for singularly perturbed problems.

Efficient numerical methods for time fractional problems have drawn the attention of several scholars over the past few decades (see \cite{Ren_Chen_AML19,Yang_et_al_AML19,Zhao_Zhang_SIAM16}). Li et al. \cite{Li_et_al_CCP18} investigated the L1-Galerkin finite element approach for solving time fractional nonlinear diffusion problems. The existence, uniqueness, and regularity for the nonlinear time-fractional diffusion problem were provided by Jin et al. \cite{Jin_et_al_SIAMJNA18}. Huang et al. \cite{Huang_etal_BIT18} established the error analysis for direct discontinuous Galerkin (DDG) finite element method for time-fractional reaction-diffusion problem.

The class of initial-boundary-value problems (IBVPs) with Caputo fractional derivative that we will be focusing on in this article is as follows:
\begin{equation}\label{ch5_eqn1}
\left\{
\begin{array}{lll}
{{}^C\mathcal{D}}^{\alpha}_{t}u(y,t)+\mathcal{L}(u(y,t))=f(y,t),\quad (y,t)\in Q:=\Omega\times (0,T],\\ [8pt]
u(y,0)=g(y), \quad\forall y\in \overline{\Omega}, \\ [8pt]
u(0,t)=u(\ell,t)=0, \quad\forall t\in (0,T],
\end{array}
\right.
\end{equation}
where $\Omega=(0,\ell),\,\mathcal{L}(u(y,t))=p(t)\left[-\dfrac{\partial}{\partial y}\left(a(y)\dfrac{\partial u}{\partial y}\right)+b(y)u\right](y,t)$ with
$ a\in \Cc^1(\overline{\Omega}),\, b\in\Cc(\overline{\Omega})$, $p\in\Cc([0,T]),\, a,p>0,\, b\geq 0$ and for $0 < \alpha < 1,\,{}^{C}\mathcal{D}^{\alpha}_tu(y,t)$ be the Caputo fractional derivative:
\begin{equation}\label{ch5_eqn2}
{}^{C}\mathcal{D}^{\alpha}_tu(y,t):=\dfrac{1}{\Gamma(1-\alpha)}\int \limits_{s=0}^t  (t-s)^{-\alpha}\frac{\partial u}{\partial s}(y,s)ds.
\end{equation}
The compatibility restrictions $g(0)=g(\ell)=0$ are likewise assumed to be valid. The presence of the unique solution to the IBVP (\ref{ch5_eqn1}) is discussed in \cite[Theorem 3.6]{MAJI2023549}.

As far as the authors knowledge goes this is the first time in the literature, NIPG method is applied to solve the time-fractional IBVPs. Section \ref{ch5_Discrete} deals with the proposed numerical scheme for the IBVP (\ref{ch5_eqn1}), whereas the $L^2$-norm stability studied in Sections \ref{ch5_sec_stab}. Section \ref{ch5_supconv} discusses the superconvergence analysis of the proposed method. By utilizing the Newton linearization approach, we expand linear problem results to semilinear problem results in Section \ref{ch5_semilinear}. Section \ref{ch5_Num_Ex} carries out numerical experiments.

\textbf{Notation.} The generic positive constant $C$ is used throughout this article, which is irrespective of step-lengths and grids, in the spatial and temporal directions. Whereas $C$ could take different values at different places. We set $p^n=p(t_n),\, g_m=g(y_m)$ and $w_m^n=w(y_m,t_n)$.

\section{The Numerical Method}\label{ch5_Discrete}

This section provides the numerical scheme to solve the IBVP (\ref{ch5_eqn1}). We further demonstrate that the associated discretized scheme has the Galerkin orthogonality property.
	
\subsection{ Temporal Semi-discretization}\label{ch5_temp_discr}

Set $t_n=T\left(\dfrac{n-1}{N}\right)^r$ for $n=1,\cdots,N+1,$ where $N$ denotes a positive integer and $r\geq 1$ is arbitrary. For $n=1,2,\cdots,N$, let $\tau_n=t_{n+1}-t_{n}$.

The Caputo fractional derivative is discretized using the standard L1-approximation as $^C\mathcal{D}^{\alpha}_{t}u(y,t_n)$:
\begin{equation}\label{ch5_L1}
^C\mathcal{D}^{\alpha}_{N}u(y,t_n):=\widehat{T}_{n,n-1}u^n(y)+\displaystyle\sum\limits_{j=2}^{n-1} [\widehat{T}_{n,j-1}-\widehat{T}_{n,j}]u^j(y)-\widehat{T}_{n,1}\,g(y),
\end{equation}
\[ \text{where}\quad
\widehat{T}_{n,j}=\dfrac{(t_n-t_j)^{1-\alpha}-(t_n-t_{j+1})^{1-\alpha}}{\tau_{j}\cdot \Gamma(2-\alpha)}, \quad 1\leq j\leq n-1,\quad n\geq 2,
\]
and $\widehat{T}_{n,j}\geq \widehat{T}_{n,{j-1}}$. From \cite[Lemma 5.2]{Stynes_et_al_SIAM2017}, the bound for the truncation error at time $t=t_n$, $\mathcal{R}^n_N(y)\displaystyle := \left({{}^C\mathcal{D}_N^{\alpha}}-{{}^C\mathcal{D}_t^{\alpha}}\right)u(y,t_n)$, is as follows:
\begin{equation}\label{ch5_Tr1}
\left|\mathcal{R}^n_N(y)\right|\leq Cn^{-{\min\{2-\alpha,r\alpha\}}}.	
\end{equation}	

We generate the following semi-discrete problem by discretizing the temporal derivative in the IBVP (\ref{ch5_eqn1}): For 
\begin{equation}\label{ch5_semi1}
\left\{
\begin{array}{lll}
\displaystyle -\dfrac{\partial}{\partial y}\left(K^n\dfrac{\partial \widehat{u}^n}{\partial y}\right)(y)+(c^n\widehat{u}^n)(y)=F^n(y),\quad y\in \Omega,\\ [9pt]
\widehat{u}^1(y)=g(y),\,y\in \Omega, \\ [9pt]
\widehat{u}^n(0)=\widehat{u}^n(\ell)=0,\quad \text{for}\quad n=1,2,\cdots,N,
\end{array}
\right.
\end{equation}
where $K^n(y)=p(t_n)a(y),\, c^n(y)=p(t_n)b(y)+\widehat{T}_{n,n-1}$, $f^n(y)=f(y,t_n)$ and
$F^n(y)=\sum\limits_{j=2}^{n-1}[\widehat{T}_{n,j} - \widehat{T}_{n,j-1}]\, \widehat{u}^j(y)+\widehat{T}_{n,1}\,\widehat{u}^1(y) + f^n(y)$.

The following conclusion can be reached by applying the same method of proof as in \cite[Lemma 4.5]{MAJI2023549}.
\begin{lemma}\label{ch5_int0}
The error resulting from temporal semi-discretization of IBVP (\ref{ch5_eqn1}) has the following bound:
\begin{equation}
\|u^n-\widehat{u}^n\|\leq CN^{-\min\{2-\alpha,r\alpha\}}, \quad n=1,2,\cdots,N+1.
\end{equation}
\end{lemma}

\subsection{The fully discrete NIPG method}

The spatial derivatives are discretized using the NIPG approach in this section. For a positive integer $M$, we discretize the spatial domain as
\[{\overline{\Omega}_M} = \left\{y_m: y_m= (m-1)h,\,\, m=1,2,\cdots,M+1 \right\},\]
where $h = \ell/M$. Let $y_{m+1/2}=(y_m+y_{m+1})/2$ for $m=1,2,\cdots,M$. Let the partition of the domain $\Omega$ represented by $ D_M = \{ \mathcal{K}_m = (y_m,y_{m+1}): m=1,2,\cdots,M \} $. We define the broken Sobolev space of order $d\geq 0$ for each element $\mathcal{K}_m \in D_M $ as
\[
H^d (\Omega,D_M) = \left\{v \in L^2(\Omega): v_{|_{\mathcal{K}_m}} \in H^d(\mathcal{K}_m),\quad \forall \mathcal{K}_m \in D_M  \right\}.
\]
The related broken Sobolev norm and semi-norm of order $d_1\leq d$ are provided by
\begin{equation*}
\norm{v}^2_{d_1,D_M} = \sum_{m=1}^{M} \norm{u}^2_{d_1,\mathcal{K}_m}, \, |u|^2_{d_1,D_M} = \sum_{m=1}^{M} |u|^2_{d_1,\mathcal{K}_m},\, \|v\|^2=\norm{v}^2_{0,D_M} =\sum\limits_{m=1}^M\,\int\limits_{\mathcal{K}_m}v^2dy,
\end{equation*}
where the normal Sobolev norm and semi-norm are defined over the domain $\mathcal{K}_m$, respectively, $\norm{\cdot}_{d_1,\mathcal{K}_m}$ and $|\cdot|_{d_1,\mathcal{K}_m}$.

For a fixed $k\geq 1$, the finite element space $V^k_{h}(\Omega)$   associated to the family $D_M$ will be defined as follows:
\begin{eqnarray*}
&&V^k_{h}(\Omega) = \left\{ v \in L^2(\Omega): v_{|_{\mathcal{K}_m}}\in \mathcal{P}^k(\mathcal{K}_m),\, \forall \mathcal{K}_m \in D_M \right\},\\
&\text{and}&\quad
V^k_{0,h}(\Omega) = \left\{ v \in V^k_h: v(0)=v(\ell)=0 \right\},
\end{eqnarray*}
where $\mathcal{P}^k(\mathcal{K}_m) $ stands for the space of polynomials of degree at most $k$ on $\mathcal{K}_m$. In order to analyze DG method, we employ the usual notations:
\[
v^{\pm}(y_m)=v(y_m\pm 0),\,\, [v(y_m)]=v^+(y_m)-v^-(y_m),\,\, \{v(y_m)\}=\dfrac{v^+(y_m)+v^-(y_m)}{2},
\]
and
\[[v(y_1)]=v(y_1),\,\,\{v(y_1)\}=v(y_1),\,\,[v(y\strut_{M+1})]=-v(y\strut_{M+1}),\,\, \{v(y\strut_{M+1})\}=v(y\strut_{M+1}).  \]

The finite element approximation for the semi-discrete problem (\ref{ch5_semi1}) is given by the following using the NIPG method for spatial variable:
\begin{equation}\label{ch5_bilnr}
\left\{
\begin{array}{ll}
\mbox{Find $u_h^n\in V^k_{0,h} (n=2,\cdots,N+1)$ so that}\\[6pt]
B(u^n_h,v_h)=\langle F^n,v_h\rangle,\quad \forall v_h\in  V^k_{0,h},\\[6pt]
\displaystyle\int\limits_{\Omega} u^1_hv_hdy=\int\limits_{\Omega} g(y)v_hdy,\quad \forall v_h\in  V^k_{0,h},
\end{array}
\right.
\end{equation}
with $B(u_h^n,v_h) = B_1(u_h^n,v_h)+B_2(u_h^n,v_h)-B_2(v_h,u_h^n)+B_3(u_h^n,v_h),$
where
\begin{equation*}
\begin{array}{lll}
\displaystyle B_1(u_h^n,v_h)=\sum\limits_{m=1}^{M}\, \int\limits_{\mathcal{K}_m}\left(K^n(y)(u_h^n)_y(v_h)_y+c^n(y) u_h^n v_h\right)dy,\\
\displaystyle B_2(u_h^n,v_h)=\sum_{m=1}^{M+1}\{ K^n(y_m) (u_h^n)_y (y_m)\}[ v_h(y_m)],\\
\displaystyle B_3(u_h^n,v_h)=\sum_{m=1}^{M+1}\sigma_m [ u_h^n(y_m)] [ v_h(y_m)],\\
\displaystyle\langle F^n,v_h\rangle = \sum_{m=1}^{M}\, \int\limits_{\mathcal{K}_m}\Big(f^n+\widehat{T}_{n,1}u_h^1+\sum\limits_{j=2}^{n-1} (\widehat{T}_{n,j}-\widehat{T}_{n,j-1})u_h^j\Big)v_hdy,
\end{array}
\end{equation*}
and the stated discontinuity-penalization parameters $\sigma_m \geq 0\,\,(m=1,2,\cdots,M+1)$ are connected to the node $y_m$.

Using the definitions of jump and average as well as integration by parts, one can now obtain the Galerkin orthogonality property for the bi-linear form $B(\cdot,\cdot)$.
\begin{lemma}\label{ch5_Galerkin_Orthognality}
Let $\widehat{u}^n$ be the solution of the semi-discrete problem (\ref{ch5_semi1}) and $u^n_h$ be the solution of the fully-discrete scheme (\ref{ch5_bilnr}). The Galerkin orthogonality property is then satisfied by the bi-linear form $B(\cdot,\cdot)$:
\[B(\widehat{u}^n-u_h^n,v_h)=0,\quad\forall v_h\in V^k_{0,h}(\Omega),\,1\leq n\leq N+1.\]
\end{lemma}

{\bf{Proof.}}
Since $\widehat{u}^n$ is the solution of (\ref{ch5_semi1}), we have $[\widehat{u}^n(y_m)]=0,\, 1\leq m\leq M+1$ and $[\widehat{u}_y^n(y_m)]=0,\, 2\leq m\leq M$. Then for all $v_h\in V^k_{0,h}(\Omega)$, we have
\begin{equation}\label{ch5_BB}
B(\widehat{u}^n,v_h)=B_1(\widehat{u}^n,v_h)+B_2(\widehat{u}^n,v_h).
\end{equation}
	
Using the definitions of jump and average as well as integration by parts, one can now obtain that
\begin{equation}\label{ch5_B1}
B_1(\widehat{u}^n,v_h)=\sum_{m=1}^{M}\, \int\limits_{\mathcal{K}_m}\Big(-(K^n(\widehat{u}^n)_y)_y+c^n\widehat{u}^n \Big)v_hdy
-\sum_{m=1}^{M+1}\{K^n(y_m)(\widehat{u}^n)_y(y_m)\}[v_h(y_m)].
\end{equation}
Thus from (\ref{ch5_BB}), (\ref{ch5_B1}) and using (\ref{ch5_semi1}), we get the following for $v_h\in V^k_{0,h}(\Omega)$:
\begin{equation}\label{ch5_orth1}
B(\widehat{u}^n,v_h)=\sum_{m=1}^{M}\, \int\limits_{\mathcal{K}_m}\Big(-(K^n(\widehat{u}^n)_y)_y+c^n\widehat{u}^n \Big)v_hdy=\langle F^n,v_h\rangle.
\end{equation}
From (\ref{ch5_bilnr}) and (\ref{ch5_orth1}), we obtained our required result. 
\hfill\eop

For the given bi-linear form $B(\cdot,\cdot)$, we define the DG energy-norm as
\[\verti{v}_{DG}^2=\beta\sum\limits_{m=1}^M\,\int\limits_{\mathcal{K}_m}(|v_y|^2+|v|^2)dy + \sum\limits_{m=1}^{M+1}\sigma_m[v(y_m)]^2,\quad v\in V^k_{0,h}(\Omega), \]
and the discrete energy-norm as
\[\verti{v}^2=\beta h \sum\limits_{m=1}^M\,\sum\limits_{j=1}^k w_j|v_y(y_{mj})|^2+\beta \|v\|^2 + \sum\limits_{m=1}^{M+1}\sigma_m[v(y_m)]^2,\quad v\in V^k_{0,h}(\Omega), \]
where $\beta=\min\{p_*a_*,b_*p_*\}$, $p_*=\min\limits_{1\leq n\leq N}p(t_n)$, $a_*=\min\limits_{y\in \Omega}a(y)$ and $y_{mj}$ are the Gaussian points in element $K_m$, $w_j>0$ are weights for the $k$-point Gaussian quadrature rule. It is obvious that the coercivity condition satisfied by the bi-linear form provided in (\ref{ch5_bilnr}) in both the energy-norm, i.e., $B(v,v)\geq \verti{v}^2\,\text{and}\,B(v,v)\geq \verti{v}_{DG}^2$.

\section{$L^2-$stability of the fully-discrete scheme}\label{ch5_sec_stab}

This section establishes the $L^2-$stability for the fully discrete L1-NIPG scheme. 

\begin{lemma}
The following bound is satisfied by the solution $u^n_h$ of the fully-discrete scheme (\ref{ch5_bilnr}):
\begin{equation}\label{ch5_stab1}	
\|u^n_h\|\leq \tau_{n-1}^{\alpha}\Gamma(2-\alpha)\Big(\|f^n\|+\widehat{T}_{n,1}\|u^1_h\|+\sum\limits_{j=2}^{n-1}(\widehat{T}_{n,j}-\widehat{T}_{n,j-1})\|u^j_h\|\Big),\, n=2,3,\cdots,N+1.
\end{equation}
\end{lemma}

{\bf{Proof.}}
Fix $n\in \{2,3,\cdots,N+1\}$. Letting $v=u^n_h$ in (\ref{ch5_bilnr}) and using the condition given for $a,b,p$, we have
\begin{eqnarray*}
\widehat{T}_{n,n-1}\|u^n_h\|^2 &\leq& |B(u^n_h,u^n_h)|=|\langle F,u^n_h\rangle| \\
&\leq& \sum\limits_{m=1}^M\,\int\limits_{\mathcal{K}_m}\Big(|f^n|+\widehat{T}_{n,1}|u^1_h|+\sum\limits_{j=2}^{n-1}(\widehat{T}_{n,j}-\widehat{T}_{n,j-1})|u^j_h|\Big)|u^n_h|dy.	 
\end{eqnarray*}
Using the Cauchy-Schwarz inequality on the right side of the aforementioned inequality results in
\begin{equation*}
\widehat{T}_{n,n-1}\|u^n_h\|^2\leq\Big(\|f^n\|+\widehat{T}_{n,1}\|u^1_h\|+\sum\limits_{j=2}^{n-1}(\widehat{T}_{n,j}-\widehat{T}_{n,j-1})\|u^j_h\|\Big)\|u^n_h\|,
\end{equation*}
which is equivalent to (\ref{ch5_stab1}). Hence the proof is done.     
\hfill\eop

For $n=2,3,\cdots,N+1$ and $j=1,2,\cdots,n-1$, we define the real numbers $\vartheta_{n,j}$ by
\begin{equation}\label{ch5_theta}
\vartheta_{n,n}=1, \quad\vartheta_{n,j}=\sum\limits_{k=1}^{n-j}\tau_{n-k}^{\beta}\vartheta_{n-k,j} \left[\widehat{T}_{n,n-k}-\widehat{T}_{n,n-k-1}\right].
\end{equation}
Observe that $\widehat{T}_{n,j}\geq \widehat{T}_{n,{j-1}}$ implies $\vartheta_{n,j}>0$ for all $n,j$.

Our investigation of stability and time-discretization error is related to the following result.
\begin{lemma}\cite[Lemma 4.4]{MAJI2023549}\label{ch5_thh}
The bound for $\vartheta_{n,j}$ is as follows:
\begin{equation}\label{ch5_thetanj}
\tau_{n-1}^{\alpha}\sum\limits_{j=1}^nj^{-\mu}\vartheta_{n,j}\,\, \leq \,\, C T^{\alpha}N^{-\mu}, \quad \mbox{for $n=2,3,\cdots,N+1$},
\end{equation}	
where $\mu$ is a parameter that satisfies the conditions $\mu\leq r\alpha$ and $C\geq \Gamma(1-\alpha)$.
\end{lemma}

Next lemma provides the stability result for the discrete scheme (\ref{ch5_bilnr}).

\begin{lemma}\label{ch5_stab2}
Let $u^n_h$ be the solution of (\ref{ch5_bilnr}). Then, we have
\begin{equation}
\|u^n_h\|\leq \|u^1_h\| + C \, T^{\alpha}\max\limits_{1\leq j\leq n} \|f^j\|, \quad \mbox{$n=2,3,\cdots,N+1$}.
\end{equation}
\end{lemma}

{\bf{Proof.}}
By replicating the proof of \cite[Lemma 4.2]{Stynes_et_al_SIAM2017} in the Lemma \ref{ch5_stab1}, we have
\begin{equation}\label{ch5_stab_2}
\displaystyle\|u_h^n\|\leq \|u^1_h\|+C\tau_{n-1}^{\alpha}\sum\limits_{j=1}^n\vartheta_{n,j}\|f^j\|
\leq \|u^1_h\|+C\tau_{n-1}^{\alpha}\max\limits_{1\leq j\leq n}\|f^j\|\sum\limits_{j=1}^n\vartheta_{n,j}.
\end{equation}	
By Lemma \ref{ch5_thh} with $\mu=0$, we obtain our required result.  
\hfill\eop

\begin{remark}
The NIPG approach is well posed, as shown by the aforementioned lemma, which states that if $f=u_h^1=0$, then one must have the trivial solution $u^n_h=0$ in (\ref{ch5_bilnr}). This also concludes the uniqueness of the discrete solution.
\end{remark}

\section{Superconvergence analysis of proposed DG method}\label{ch5_supconv}
In this section, we will study the superconvergence analysis of the fully discrete scheme using discrete energy-norm rather than using the DG energy-norm. Let us first consider the Lobatto points $\{-1=z_0^{(k)},z_1^{(k)},\cdots,z_{k-1}^{(k)},z_k^{(k)}=1\}$ on $[-1,1]$, which are the $(k+1)$ zeros of the polynomial $\Phi_{k+1}(z)=\dfrac{d^{k-1}}{dz^{k-1}}(z^2-1)^k$.
Let $\Pi v\in \mathcal{P}^k([-1,1])$ denotes the Lagrange interpolation of a continuous function $v$ defined on $[-1,1]$ by 
using the Lobatto nodes $\{z_j^{(k)}\}_{j=0}^k$. Then we have
\begin{eqnarray}
(\Pi v-v)(z)&=&\dfrac{v^{(k+1)}(z)}{(2k)!}\Phi_{k+1}(z)+O(1)\int\limits_{-1}^{1}|v^{(k+2)}(z)|dz, \label{ch5_interp1}\\ [8pt]
(\Pi v-v)'(z)&=&\dfrac{v^{(k+1)}(z)}{(2k)!}\Psi_{k}(z)+O(1)\int\limits_{-1}^{1}|v^{(k+2)}(z)|dz, \label{ch5_interp2}
\end{eqnarray}
where $\Psi_{k}(z)=\dfrac{d^{k}}{dz^{k}}(z^2-1)^k$ is a multiple of the $k$th Legendre polynomial on $[-1,1]$ and Gauss points are the zeros of $\Psi_{k}(z)$.

For $m=1,2,\cdots,M$, let $y_{m,j}^{(k)}=\dfrac{y_m+y_{m+1}}{2}+\dfrac{y_{m+1}-y_m}{2}z_j^{(k)},\, j=0,1,\cdots,k$ are the nodes in $[y_m,y_{m+1}]$ corresponding to $z_j^{(k)}$. If $\pi_h u $ denotes the piecewise Lagrange interpolation of $\widehat{u}^n$ using $\{y_{m,j}^{(k)}\}_{j=0}^k$ as node points on $[y_m,y_{m+1}]$, then from equations (\ref{ch5_interp1}) and (\ref{ch5_interp2}), we have the following for $y\in [y_m,y_{m+1}]$ and $n=1,2,\cdots,N+1$:
\begin{eqnarray}
\displaystyle(\pi_h\widehat{u}^n-\widehat{u}^n)(y)&=&\dfrac{(\widehat{u}^n)^{(k+1)}(y)}{(2k)!}\widehat{\Phi}_{m,k+1}(y)+O(h^{k+1})\int\limits_{\mathcal{K}_m}\left|(\widehat{u}^n)^{(k+2)}(y)\right|dy, \label{ch5_interp3}\\ [8pt]
(\pi_h\widehat{u}^n-\widehat{u}^n)'(y)&=&\dfrac{(\widehat{u}^n)^{(k+1)}(y)}{(2k)!}\widehat{\Psi}_{m,k}(y)+O(h^k)\int\limits_{\mathcal{K}_m}\left|(\widehat{u}^n)^{(k+2)}(y)\right|dy, \label{ch5_interp4}
\end{eqnarray}
where $\widehat{\Phi}_{m,k+1}(y)=\dfrac{d^{k-1}}{dy^{k-1}}\left((y-y_{m+1/2})^2-\dfrac{h^2}{4}\right)^k$ and $\widehat{\Psi}_{m,k}(y)=\widehat{\Phi}'_{m,k+1}(y)$.

Therefore, the equation (\ref{ch5_interp4}) at the Gauss points $y_{mj}$ becomes
\begin{equation}\label{ch5_interp5}
(\pi_h\widehat{u}^n-\widehat{u}^n)'(y_{mj})=O(h^k)\int\limits_{\mathcal{K}_m}|(\widehat{u}^n)^{(k+2)}(y)|dy.
\end{equation}

We have the following interpolation error estimate for any $\widehat{u}^n\in H^{k+1}(\Omega,D_M)$ ($n=1,2,\cdots,N+1$):
\begin{equation}\label{ch5_result_int1}
\|\widehat{u}^n-\pi_h\widehat{u}^n\|_{d_1,D_M} \leq Ch^{k+1-d_1}|\widehat{u}^n|_{k+1,D_M},\quad 0\leq d_1\leq k+1.
\end{equation}

Let $u^n$ and $u_h^n$ respectively, denote  the solutions at $t=t_n$ of the equations (\ref{ch5_eqn1}) and (\ref{ch5_bilnr}), then we define the error $e^n$ as
\[
e^n=u^n-u^n_h=(u^n-\widehat{u}^n)+(\widehat{u}^n-\pi_h\widehat{u}^n)+(\pi_h\widehat{u}^n-u^n_h),\, n=1,2,\cdots,N+1. \]
Let's use the notation $\eta^n=\widehat{u}^n-\pi_h\widehat{u}^n$ for the interpolation error and $\xi^n=\pi_h\widehat{u}^n-u^n_h$ for the spatial discretization error.

The interpolation error in the DG energy-norm will be investigated in the following theorem.
\begin{theorem}\label{ch5_int_err1}
The interpolation error can be bounded by
\begin{equation*}
\verti{\eta^n}_{DG}\leq Ch^k, \quad n=1,2,\cdots,N+1.
\end{equation*}
\end{theorem}

{\bf Proof.}
Now, since $\eta^n(y_m)=0$ for $m=1,2,\cdots,M+1,$ which implies that $[\eta^n(y_m)]=0$. Then for each $n=1,2,\cdots,N+1$, we have
\begin{equation}\label{ch5_eta_bound}
\verti{\eta^n}_{DG}^2=\beta\sum\limits_{m=1}^M\,\int\limits_{\mathcal{K}_m}(|\eta^n_y|^2+|\eta^n|^2)dy=\beta\Big(|\eta^n|_{1,D_M}^2+\|\eta^n\|^2\Big)\leq C h^{2k},
\end{equation}
where we have used the estimate given in (\ref{ch5_result_int1}). Hence the proof is done. \hfill \eop

The interpolation error in discrete energy-norm is a topic we'll explore in the following theorem.
\begin{theorem}\label{ch5_int_err}
The interpolation error can be bounded by
\begin{equation*}
\verti{\eta^n}\leq C h^{k+1}, \quad n=1,2,\cdots,N+1.
\end{equation*}
\end{theorem}

{\bf{Proof.}}
Now, since $\eta^n(y_m)=0$ for $m=1,2,\cdots,M+1,$ which implies that $[\eta^n(y_m)]=0$. Then for each $n=1,2,\cdots,N+1$, we have
\begin{equation}\label{ch5_eta_energy}
\verti{\eta^n}^2=\beta \sum\limits_{m=1}^M\,h \sum\limits_{j=1}^k w_j |(\eta^n)'(y_{mj})|^2+\beta\|\eta^n\|^2.
\end{equation}
By using equations (\ref{ch5_interp5}) and (\ref{ch5_result_int1}) in the above equation (\ref{ch5_eta_energy}) and then using H\"{o}lder's inequality, we have the following
\begin{eqnarray*}
\verti{\eta^n}^2&\leq&C\sum\limits_{m=1}^M h^{2k+1} \left(\,\,\int\limits_{\mathcal{K}_m}|(\widehat{u}^n)^{(k+2)}(y)|dy\right)^2+Ch^{2k+2}\left(\sum\limits_{m=1}^M\,\int\limits_{\mathcal{K}_m}|(\widehat{u}^n)^{(k+1)}(y)|dy\right)^2\\ [8pt]
&\leq& C h^{2k+2}\,|\widehat{u}^n|_{k+2,D_M}^2+C h^{2k+3}\,|\widehat{u}^n|_{k+1,D_M}^2\\
&\leq& Ch^{2k+2}.
\end{eqnarray*}
Hence the proof is done.
\hfill\eop

In the next we will find out the bound for the discretization error in discrete energy-norm.
\begin{theorem}\label{ch5_dis_err}
The following bound is satisfied by the discretization error:
\begin{equation*}
\verti{\xi^n}\leq Ch^{k+1}, \quad  n=1,2,\cdots,N+1.
\end{equation*}
\end{theorem}

{\bf{Proof.}}
Using Lemma \ref{ch5_Galerkin_Orthognality} and the coercivity condition, we have
\begin{equation}\label{ch5_xi_1}
\verti{\xi^n}^2\leq B(\xi^n,\xi^n)=B(\widehat{u}^n-u_h^n,\xi^n)-B(\eta^n,\xi^n)=-B(\eta^n,\xi^n).
\end{equation}
As $[\eta^n(y_m)]=0$, we can estimate $B_2(\xi^n,\eta^n)=B_3(\eta^n,\xi^n)=0.$ So $B(\eta^n,\xi^n)=B_1(\eta^n,\xi^n)+B_2(\eta^n,\xi^n)$.
Now by using H\"{o}lder's inequality, $B_1(\eta^n,\xi^n)$ can be written as
\begin{eqnarray}\label{ch5_B1_1}
B_1(\eta^n,\xi^n)&=&\sum_{m=1}^{M}\, \int\limits_{\mathcal{K}_m}\Big(K^n(y)\eta^n_y\xi^n_y+c^n(y)\eta^n\xi^n\Big)dy\nonumber\\
&\leq& C\Big(\|\eta^n\|_{1,D_M}\|\xi^n\|_{1,D_M}+\|\eta^n\|\|\xi^n\|\Big)
\leq Ch^{k+1}\verti{\xi^n}.
\end{eqnarray}

The remaining term $B_2(\eta^n,\xi^n)$ can be bounded by
\begin{equation}\label{ch5_B2_1}
B_2(\eta^n,\xi^n)=\sum_{m=1}^{M+1}\{ K^n(y_m) \eta^n_y (y_m)\}[ \xi^n(y_m)]\leq C h^{k+1}\verti{\xi^n}.
\end{equation}

Therefore, by combining (\ref{ch5_xi_1}) and (\ref{ch5_B1_1}), we have
\begin{equation}\label{ch5_xi_2}
\verti{\xi^n}^2\leq|B(\eta^n,\xi^n)|\leq |B_1(\eta^n,\xi^n)|+|B_2(\eta^n,\xi^n)|\leq C h^{k+1}\verti{\xi^n}.
\end{equation}
Hence, we can obtain the required result.   
\hfill\eop

The error bound resulting from the IBVP (\ref{ch5_eqn1})'s fully discretization is shown in the following theorem.
\begin{theorem}\label{ch5_full_err}
Let $u^n$ be the solution of the IBVP (\ref{ch5_eqn1}) and $u_h^n$ be the solution of the fully discrete scheme (\ref{ch5_bilnr}) at $t=t_n$, then the error estimate is as follows:
\[
\max\limits_{1\leq n\leq N+1}\verti{u^n-u^n_h}\leq C \left(h^{k+1}+N^{-\min\{2-\alpha,r\alpha\}}\right).
\]
\end{theorem}

{\bf{Proof.}}
We know that for each $n=1,2,\cdots,N+1$,
\[\verti{u^n-u^n_h}\leq \verti{u^n-\widehat{u}^n}+\verti{\eta^n}+\verti{\xi^n}.\]
Now by combining the results given in Lemma \ref{ch5_int0}, Theorem \ref{ch5_int_err} and Theorem \ref{ch5_dis_err} and we derive the necessary estimate by considering maximum over $n$.
\hfill\eop

\begin{remark}
In a similar way of Theorem \ref{ch5_dis_err} and Theorem \ref{ch5_full_err}, one can establish that $\verti{\xi^n }_{DG}\leq Ch^k$ and $\max\limits_{1\leq n\leq N+1}\verti{u^n-u^n_h }_{DG}\leq Ch^k$. As a result, we have a one-order increase in the order of convergence when utilizing the discrete energy-norm rather than the DG energy-norm.
\end{remark}

\section{Semi-linear time-fractional IBVPs}\label{ch5_semilinear}

We consider the following class of semi-linear time-fractional IBVPs in this section:
\begin{equation}\label{ch5_smlnr_eqn1}
	\left\{
	\begin{array}{lll}
		 {{}^C\mathcal{D}}^{\alpha}_{t}u(y,t)+\mathcal{L}(u(y,t))=f(y,t),\quad (y,t)\in Q:=\Omega\times (0,T],\\ [8pt]
		u(y,0)=g(y), \quad\forall y\in \overline{\Omega}, \\ [8pt]
		u(0,t)=u(\ell,t)=0, \quad\forall t\in (0,T],
	\end{array}
	\right.
\end{equation}
where $\Omega=(0,\ell)$ and $\mathcal{L}(u(y,t))=p(t)\left[-\dfrac{\partial}{\partial y}\left(a(y)\dfrac{\partial u}{\partial y}(y,t)\right)+b(y,u)\right]$ with
\begin{equation}\label{ch5_smlnr_con_coef}
	\left\{
	\begin{array}{ll}		
		a\in \Cc^1(\overline{\Omega}),\,\, b\in\Cc^1(\Omega\times \mathbb{R}),\,\, a(y)>0,\,\,b_{u}\geq 0, \quad\forall y\in\overline{\Omega},\\ [7pt]
		p\in\Cc([0,T]),\,\, p(t)> 0, \quad \forall t\in (0,T].
	\end{array}\right.
\end{equation}	
We further assume that the problem (\ref{ch5_smlnr_eqn1}) admits an unique solution $u(y,t)$ and that the compatibility conditions $g(0)=g(\ell)=0$ are true for suitably smooth functions $b(y,u),f(y,t)$ and $g(y)$.

In order to arrive at the problem's numerical solution, we first utilize the Newton linearization process to convert the semi-linear problem (\ref{ch5_smlnr_eqn1}) into a sequence of linear problems and then we use the numerical approach described in the above section. As a result, starting with the first guess $u^{(0)}$, we define the solution $u^{(q+1)},\,q \geq 0$ of the following linear time-fractional IBVP:
\begin{equation}\label{ch5_lnr}
	\left\{
	\begin{array}{lll}
		 {}^{C}\mathcal{D}_{t}^{\alpha}u^{(q+1)}(y,t)+\mathcal{L}^q(u^{(q+1)}(y,t))
		=F^{(q)}(y,t),\quad\forall (y,t) \in Q, \\ [8pt]
		u^{(q+1)}(y,0)=g(y), \quad\forall y\in \overline{\Omega}, \\ [8pt]
		u^{(q+1)}(0,t)=u^{(q+1)}(\ell,t)=0, \quad\forall t\in (0,T],
	\end{array}
	\right.
\end{equation}
where
\[
\begin{array}{ll}
\mathcal{L}^q(u^{(q+1)}(y,t))=p(t)\left[-\dfrac{\partial}{\partial y}\left(a(y)\dfrac{\partial u^{(q+1)}}{\partial y}\right) +b_{u}(y,u^{(q)})u^{(q+1)}\right](y,t),\\[12pt]
F^{(q)}(y,t)=f(y,t)-p(t)\left[b(y,u^{(q)})-b_{u}(y,u^{(q)})u^{(q)}(y,t)\right].
\end{array}
\]

We choose the convergence condition for $\{u^{q}\}_{q \geq 0}$ as:
\[
\displaystyle\max\limits_{(y_m,t_n)}|u^{(q+1)}(y_m,t_n)-u^{(q)}(y_m,t_n)| \,\, \leq \,\, \mbox{Tol},\quad q\geq 0.
\]
We used $\mbox{Tol} = 10^{-7}$ for computational purposes. The semi linear problem's numerical solutions are presented in the following section (\ref{ch5_smlnr_eqn1}).

\goodbreak\noindent
\section{Numerical tests}\label{ch5_Num_Ex}

To demonstrate the effectiveness and precision of the suggested numerical methodology, we conduct numerical experiments in this section. We consider two test examples, one with constant coefficients and another variable coefficients. Tables and figures are used to present the numerical results.

The $L^{\infty}(L^2)$ error and the order of convergence can be calculated using
\[E_{M,N}=\max\limits_{1\leq n\leq N+1} \|u^n-u^n_h\|, \quad \mbox{and} \quad q\strut_{M,N}=\log_2\left(\dfrac{E_{M,N}}{E_{2M,2N}}\right).\]
For the numerical studies, we utilized $r = ({2-\alpha})/{\alpha}$.

\begin{example}\label{ch5_ex2}
Consider the constant coefficient example over $(0,\pi)\times (0,1]$:
\begin{equation}\label{ch5_ex2_eq1}
\left\{
\begin{array}{ll}
{^C\mathcal{D}}^{\alpha}_{t}u(y,t)-\dfrac{\partial^2 u}{\partial y^2}(y,t)+u(y,t)=f(y,t),\\ [8pt]
u(y,0)=0, \quad y\in [0,\pi], \\ [8pt]
u(0,t)=u(\pi,t)=0, \quad t\in (0,1].
\end{array}
\right.
\end{equation}
\end{example}

In the IBVP (\ref{ch5_ex2_eq1}), the function $f(y,t)$ is computed such that $u(y,t)=(t^{\alpha}+t^3)\sin y$.

\begin{figure}[htb]
\centerline{
\begin{tabular}{cc}
\resizebox*{9cm}{!}{\includegraphics{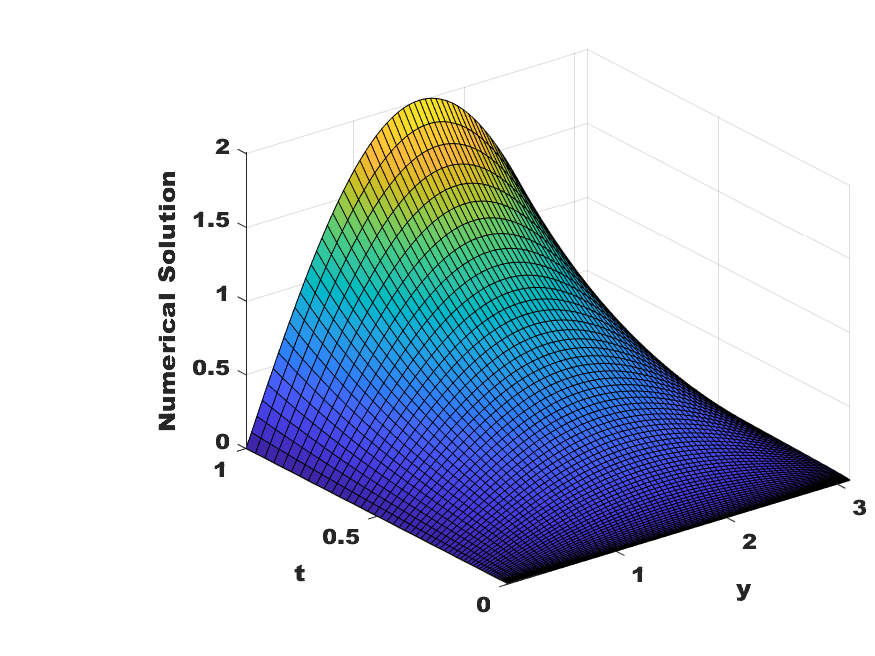}}&
\resizebox*{9cm}{!}{\includegraphics{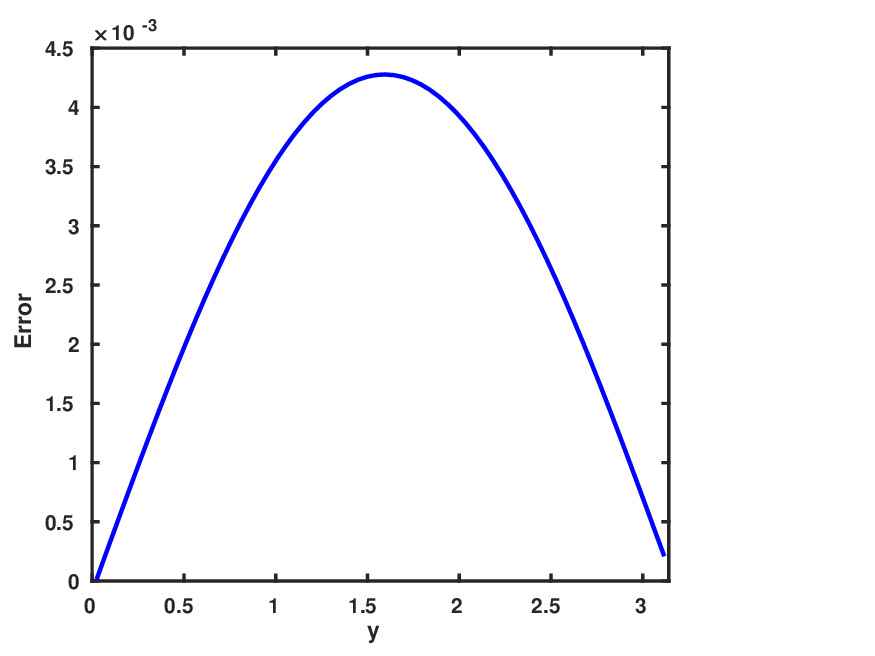}}\\
{\it (a) Numerical solution.} & {\it (b) Error at $t=1$.}
\end{tabular}
}
\caption{{\label{ch5_ex2_figure1} \it Solution surface and error curve at $t=1$ for Example \ref{ch5_ex2}.}}
\end{figure}

The numerical solution's surface plot and the error curve at $t=1$ of Example \ref{ch5_ex2} for $\alpha=0.6$ with $M=N=64$ are displayed in Figure \ref{ch5_ex2_figure1}.

The $L^2(\Omega)$, $L^\infty(\Omega)$, $\verti{\cdot}_{DG}$ and $\verti{ \cdot }$ errors and corresponding rate of convergence at $t=1$ for $\alpha=0.4,0.6$ for Example \ref{ch5_ex2} are provided in Table \ref{ch5_ex2_tab1}, where we observe that superconvergence of error occurring in $\verti{\cdot}$ over $\verti{\cdot}_{DG}$.

\begin{table}[h!]
\caption{\label{ch5_ex2_tab1}  \it{Errors and order of convergence at $t=1$ for Example \ref{ch5_ex2} with $N=\left[M^\frac{2}{2-\alpha}\right]$.}} \vspace{0.1cm}
{\centering
\begin{tabular}{c|c|ccccc}
\hline
$\alpha\downarrow$& $M\rightarrow$&$20$& $40$&$80$&$160$&$320$   \\ 
\hline 
&$\|u-u_h\|_{L^2(\Omega)}$ &1.3522e-02&	3.5314e-03&	9.1783e-04&	2.3478e-04&	5.9674e-05\\ [4pt]
&Order &1.9370 &  1.9439 & 1.9669 & 1.9761&$-$\\ \cline{2-7} 
&$\|u-u_h\|_{L^\infty(\Omega)}$&1.0732e-02&	2.8109e-03&	7.3151e-04&	1.8722e-04&	4.7601e-05\\ [4pt]
&Order& 1.9328& 1.9421& 1.9661& 1.9757&$-$\\ \cline{2-7}
0.4&$\verti{u-u_h}_{DG}$& 1.9758e-01&9.8536e-02&	4.9231e-02&	2.4611e-02&	1.2305e-02\\ [4pt]
&Order& 1.0037 & 1.0011& 1.0003& 1.0001&$-$\\ \cline{2-7}
&$\verti{u-u_h}$& 1.5666e-02&	4.1751e-03&	1.1024e-03&	2.8450e-04&	7.2755e-05\\ [4pt]
&Order&  1.9077  &  1.9212  &  1.9541  &  1.9673&$-$\\
\hline
&$\|u-u_h\|_{L^2(\Omega)}$  &1.0765e-02&2.7140e-03&	6.8132e-04&	1.7079e-04& 4.2758e-05\\ [4pt]
&Order &1.9879& 1.9940& 1.9961& 1.9980&$-$\\ \cline{2-7}
&$\|u-u_h\|_{L^\infty(\Omega)}$&8.5329e-03 &	2.1588e-03 &	5.4281e-04 &	1.3617e-04 &	3.4104e-05\\ [4pt]
&Order&1.9828 &    1.9917  &   1.9950   &  1.9974&$-$\\ \cline{2-7}
0.6&$\verti{u-u_h}_{DG}$& 1.9725e-01 &	9.8484e-02 &	4.9224e-02 &	2.4610e-02 &	1.2304e-02\\ [4pt]
&Order&1.0020   &  1.0005   &  1.0001   &  1.0000&$-$\\ \cline{2-7}
&$\verti{u-u_h}$& 1.1500e-02 &	2.9193e-03 &	7.3534e-04 &	1.8468e-04 &	4.6279e-05\\ [4pt]
&Order& 1.9779   &  1.9892   &  1.9934   &  1.9966&$-$\\
\hline
\end{tabular}
\par}
\end{table}

\begin{example}\label{ch5_ex1}
Consider the variable coefficient example over $(0,\pi)\times (0,1]$:
\begin{equation}\label{ch5_ex1_eq1}
\left\{
\begin{array}{lll}
{^C\mathcal{D}}^{\alpha}_{t}u(y,t)+(t^2+1)\left[-\dfrac{\partial}{\partial y}\left((y+1)\dfrac{\partial u}{\partial y}\right)+u\right](y,t)=f(y,t),\\ [8pt]
u(y,0)=0, \quad y\in [0,\pi], \\ [8pt]
u(0,t)=u(\pi,t)=0, \quad t\in (0,1].
\end{array}
\right.
\end{equation}
\end{example}

In the IBVP (\ref{ch5_ex1_eq1}), the exact solution is $u(y,t)=t^\alpha y\sin y$, which will be utilized to determine the value of the function $f(y,t)$.

\begin{figure}[htb]
\centerline{
\begin{tabular}{cc}
\resizebox*{9cm}{!}{\includegraphics{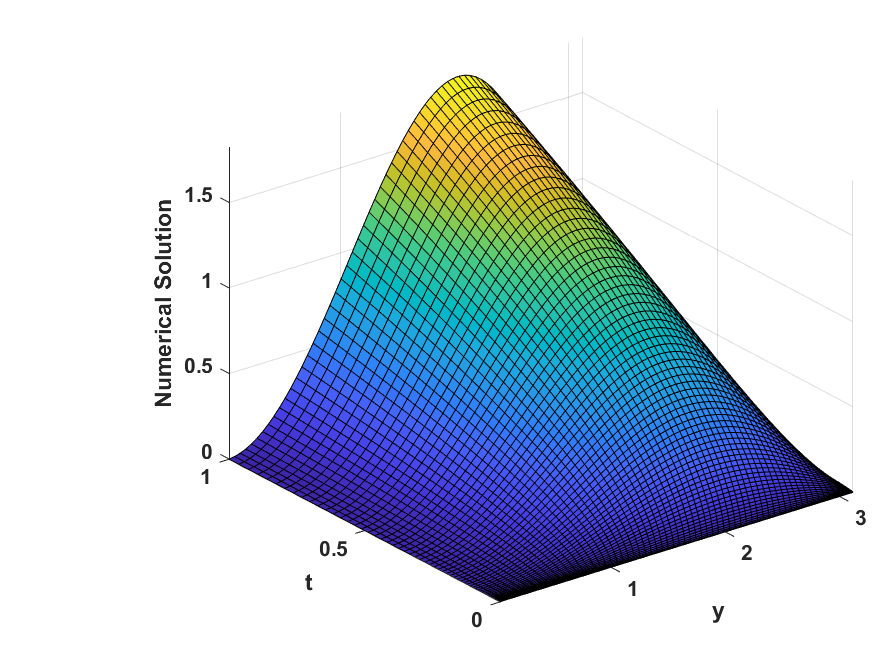}}&
\resizebox*{9cm}{!}{\includegraphics{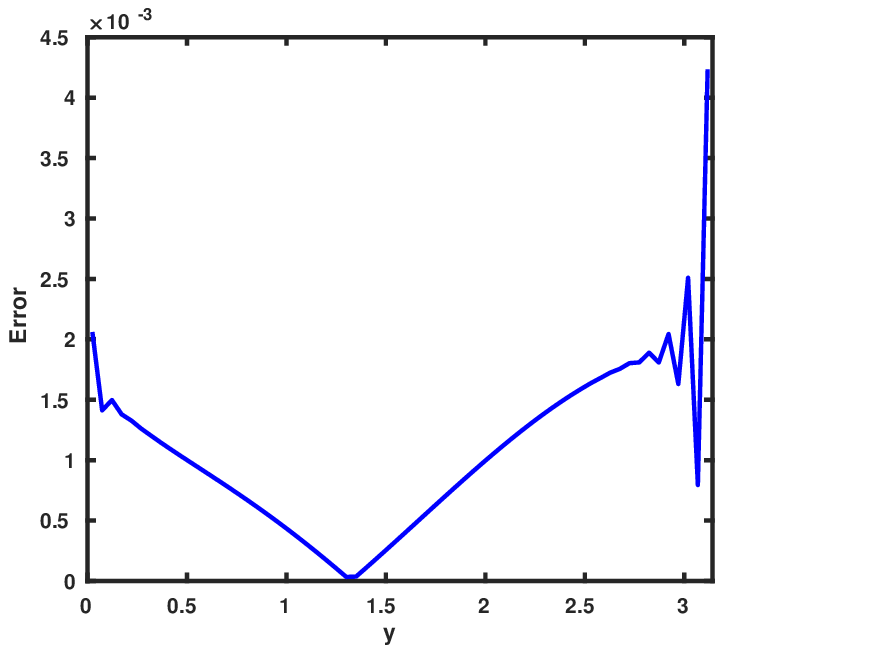}}\\
{\it (a) Numerical solution.} & {\it (b) Error at $t=1$.}
\end{tabular}
}
\caption{{\label{ch5_ex1_figure1} \it Solution surface and error curve at $t=1$ for Example \ref{ch5_ex1}.}}
\end{figure}

The numerical results for Example \ref{ch5_ex1} are given in Figure \ref{ch5_ex1_figure1}, which depicts the numerical solution's surface plot, and the error curve at $t=1$ for $\alpha=0.6$ with $M=N=64$. Also, the $L^2(\Omega)$, $L^\infty(\Omega)$, $\verti{\cdot}_{DG}$ and $\verti{\cdot}$ errors and corresponding rate of convergence at $t=1$ for $\alpha=0.4,0.6$ for Example \ref{ch5_ex1} are provided in Table \ref{ch5_ex1_tab1}, where we observe that superconvergence of error occurring in $\verti{\cdot}$ over $\verti{\cdot}_{DG}$.

\begin{table}[h!]
\caption{\label{ch5_ex1_tab1}  \it{Errors and order of convergence at $t=1$ for Example \ref{ch5_ex1} with $N=\left[M^\frac{2}{2-\alpha}\right]$.}} \vspace{0.1cm}
{\centering
\begin{tabular}{c|c|ccccc}
\hline
$\alpha\downarrow$& $M\rightarrow$&$20$& $40$&$80$&$160$&$320$   \\ 
\hline 
&$\|u-u_h\|_{L^2(\Omega)}$ &6.0297e-02&	1.3814e-02&	3.2707e-03&	7.9315e-04&	1.9511e-04\\ [4pt]
&Order &2.1259  &  2.0785   & 2.0439 &   2.0233&$-$\\ \cline{2-7} 
&$\|u-u_h\|_{L^\infty(\Omega)}$&1.0348e-01&	2.6196e-02&	6.5745e-03&	1.6458e-03&	4.1167e-04\\ [4pt]
&Order& 1.9819 & 1.9944 &   1.9981   & 1.9993&$-$\\ \cline{2-7}
0.4&$\verti{u-u_h}_{DG}$& 3.2440e-01&	1.5317e-01&	7.4670e-02&	3.6896e-02&	1.8342e-02\\ [4pt]
&Order&1.0827  &  1.0365   & 1.0171  &  1.0083&$-$\\ \cline{2-7}
&$\verti{u-u_h}$& 1.6262e-01&	5.5796e-02&	1.9477e-02&	6.8469e-03&	2.4143e-03\\ [4pt]
&Order& 1.5433  &  1.5184  &  1.5082 &   1.5038&$-$\\
\hline
&$\|u-u_h\|_{L^2(\Omega)}$ &6.0228e-02&	1.3801e-02&	3.2678e-03&	7.9247e-04&	1.9495e-04\\ [4pt]
&Order &2.1257  &  2.0784  &  2.0439  &  2.0233&$-$\\ \cline{2-7}
&$\|u-u_h\|_{L^\infty(\Omega)}$&1.0342e-01&	2.6192e-02&	6.5742e-03&	1.6458e-03&	4.1167e-04\\ [4pt]
&Order& 1.9814  &  1.9942  &  1.9980   & 1.9992&$-$\\ \cline{2-7}
0.6&$\verti{u-u_h}_{DG}$& 3.2440e-01&	1.5317e-01&	7.4671e-02&	3.6896e-02&	1.8342e-02\\ [4pt]
&Order& 1.0826  &  1.0365  &  1.0171   & 1.0083&$-$\\ \cline{2-7}
&$\verti{u-u_h}$&1.6250e-01&	5.5780e-02&	1.9474e-02&	6.8466e-03&	2.4143e-03\\ [4pt]
&Order& 1.5426 &   1.5182  &  1.5081&    1.5038&$-$\\
\hline
\end{tabular}
\par}
\end{table}

\begin{example}\label{ch5_ex3}
Consider the semilinear problem:
\begin{equation}\label{ch5_ex3_eq1}
\left\{
\begin{array}{ll}
{^C\mathcal{D}}^{\alpha}_{t}u(y,t)-\dfrac{\partial^2 u}{\partial y^2}-\exp(-u(y,t))=f(y,t),\, (y,t)\in (0,1)\times (0,1],\\ [8pt]
u(y,0)=y^2-y, \quad y\in [0,1], \\ [8pt]
u(0,t)=u(1,t)=0, \quad t\in (0,1].
\end{array}
\right.
\end{equation}
\end{example}

In the problem (\ref{ch5_ex3_eq1}), the evaluation of $f(y,t)$ results in the exact solution, $u(y,t)=(t^\alpha+t^3+1)(y^2-y)$. In order to numerically solve this problem, we first linearized the semi-linear problem (\ref{ch5_ex3_eq1}) using Newton's approach described as in equation (\ref{ch5_lnr}), which results in the following linear problems:
\begin{equation}\label{ch5_ex3_lnr}
\left\{
\begin{array}{lll}
({}^{C}\mathcal{D}_{t}^{\alpha}u^{(q+1)}+\mathcal{L}^q(u^{(q+1)}))(y,t)=F^{(q)}(y,t),\quad\forall\, (y,t) \in (0,1)\times (0,1], \\ [8pt]
u^{(q+1)}(y,0)=0, \quad y\in [0,1], \\ [8pt]
u^{(q+1)}(0,t)=u^{(q+1)}(1,t)=0, \quad t\in [0,1],
\end{array}
\right.
\end{equation}
where for $q\geq 0$,
\[
\begin{array}{ll}
\mathcal{L}^q(u^{(q+1)})=-\dfrac{\partial^2 u^{(q+1)}}{\partial y^2}+\exp(-u^{(q)})u^{(q+1)},\\[12pt]
F^{(q)}=f+\exp(-u^{(q)})(1+u^{(q)}).
\end{array}
\]

The approximate solution to the problem (\ref{ch5_ex3}) is then obtained by using the suggested numerical technique to (\ref{ch5_ex3_lnr}). We utilize the initial supposition $u^{(0)}(y,t)=0,\,(y,t) \in (0,1)\times (0,1]$ for computation.

\begin{figure}[htb]
\centerline{
\begin{tabular}{cc}
\resizebox*{9cm}{!}{\includegraphics{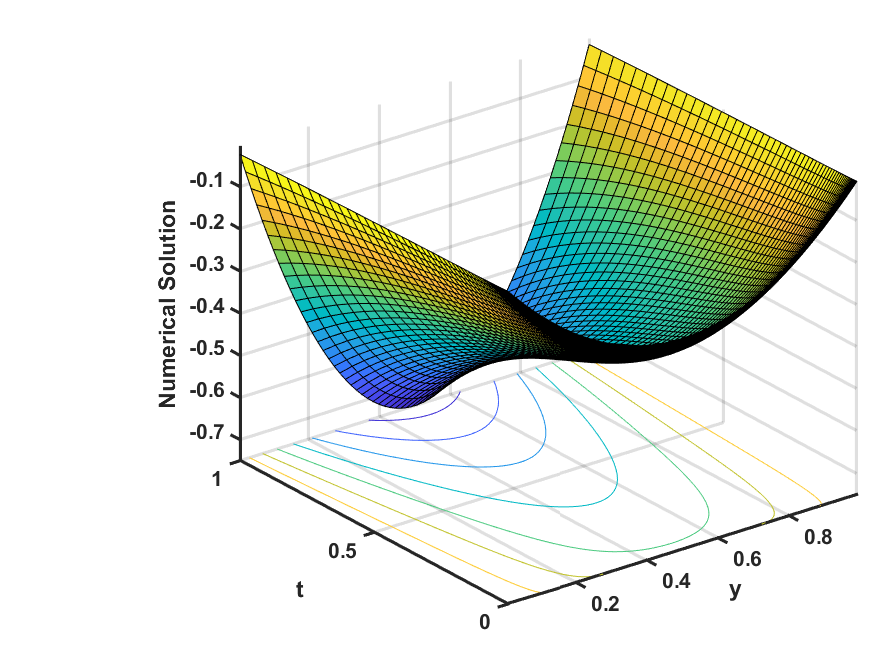}}&
\resizebox*{9cm}{!}{\includegraphics{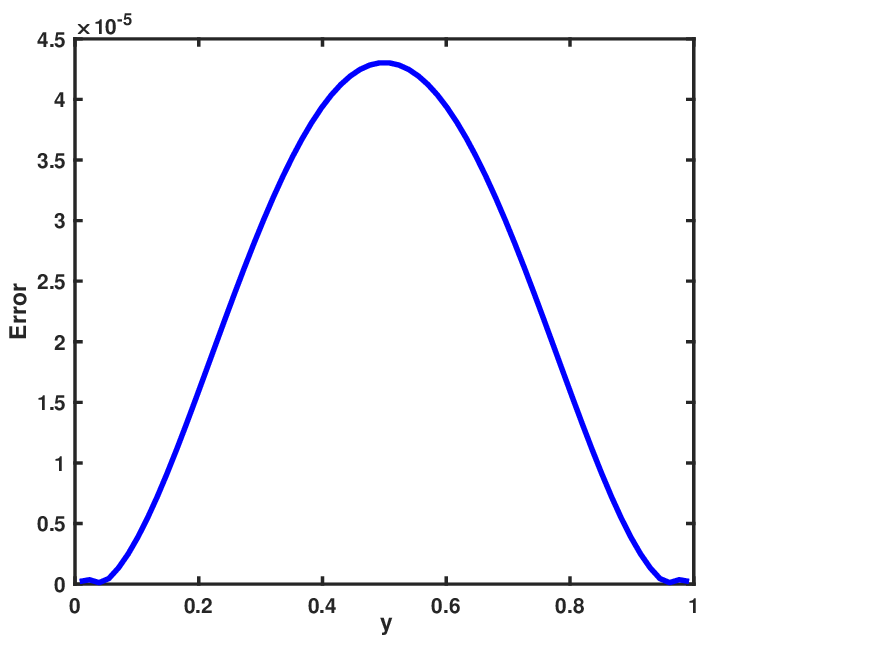}}\\
{\it (a) Numerical solution.} & {\it (b) Error at $t=1$.}
\end{tabular}
}
\caption{{\label{ch5_ex3_figure1} \it Solution surface and error curve at $t=1$ for Example \ref{ch5_ex3}.}}
\end{figure}

The numerical solution surface and error curve at $t=1$ are shown in Figure \ref{ch5_ex3_figure1} with $M=N=64$ and $\alpha=0.5$. The numerical results for the suggested approach for Example \ref{ch5_ex3} are presented in Table \ref{ch5_ex3_tab1}.

\begin{table}[h!]
\caption{\label{ch5_ex3_tab1}  \it{The $L^{\infty}(L^2)$ error and the order of convergence for Example \ref{ch5_ex3} with $k=2$.}} \vspace{0.1cm}
{\centering
\begin{tabular}{c|c|ccccc}
\hline
$\alpha\downarrow$& $M=N\rightarrow$&$16$& $32$&$64$&$128$&$256$   \\
\hline
0.4	&$E_{M,N}$ & 1.3177e-03&5.4129e-04&	2.0350e-04&	 7.3971e-05&2.6366e-05\\ [4pt]
& $q\strut_{M,N}$&1.2835 & 1.4114 & 1.4600 & 1.4883& $-$  \\
\hline
0.6	&$E_{M,N}$ &1.9151e-03&	8.8080e-04&	3.7954e-04&	 1.5574e-04&6.2081e-05 \\ [4pt]
&  $q\strut_{M,N}$&  1.1206 & 1.2146 & 1.2851 & 1.3269 &$-$ \\
\hline
0.8&$E_{M,N}$ & 2.4319e-03&	1.2029e-03&	5.7275e-04&	 2.6515e-04&1.2052e-04\\ [4pt]
&  $q\strut_{M,N}$&1.0155 & 1.0706 & 1.1111 & 1.1375& $-$\\
\hline
\end{tabular}
\par}
\end{table}

\goodbreak\noindent
\section*{Conclusions}
This article proposed and analyzed the superconvergence of the NIPG method for the non-autonomous TF-ADE (\ref{ch5_eqn1}). By discretizing the spatial derivatives using the NIPG approach and the temporal derivative using the L1-scheme, we were able to produce a fully discrete scheme. The suggested method is convergent with respect to the discrete energy-norm of order $(k+1)$ in spatial direction and of order $\min\{{2-\alpha,r\alpha}\}$ in temporal direction, where $k$ represents the degree of the piecewise polynomial in the finite element space. Additionally, after linearizing semilinear problems using the Newton linearization approach, we solved them using the suggested method. We validated the theoretical results through various numerical experiments.

\goodbreak\noindent
\section*{Acknowledgments}
The first author would like to acknowledge IIT Guwahati for their assistance with the fellowship and amenities throughout his research.

\section*{References}	
\goodbreak\noindent

\end{document}